# Wavelet block thresholding for samples with random design: a minimax approach under the $L^p$ risk

**Christophe Chesneau**

*Laboratoire de Probabilités et Modèles Aléatoires,*
*CNRS-UMR 7599, Université Paris VI, UFR de Mathématiques,*
*175 rue de Chevaleret F-75013 Paris, France.*
*e-mail:* chesneau@proba.jussieu.fr
*url:* http://www.chesneau-stat.com

**Abstract:** We consider the regression model with (known) random design. We investigate the minimax performances of an adaptive wavelet block thresholding estimator under the $\mathbb{L}^p$ risk with $p \geq 2$ over Besov balls. We prove that it is near optimal and that it achieves better rates of convergence than the conventional term-by-term estimators (hard, soft,...).

## 1. Motivations

In recent years, wavelet thresholding procedures have been widely applied to the field of nonparametric function estimation. They excel in the areas of spatial adaptivity, computational efficiency and asymptotic optimality. Among the various thresholding techniques studied in the literature, there are the term-by-term thresholding (hard, soft, ...) initially developed by Donoho and Johnstone (1995) and the block thresholding (global, BlockShrink, ...) introduced by Kerkyacharian et al. (1996) and Hall et al. (1999).

Several recent works demonstrated that the block thresholding methods can enjoy better theoretical properties than the conventional term-by-term thresholding methods. This superiority has been proved for various statistical models via the minimax approach under the $\mathbb{L}^2$ risk. See, for instance, Cai (1999) and Cavalier and Tsybakov (2001) for the Gaussian sequence model, Cai and Chicken (2005) for the density estimation, Chicken (2003) for the regression model with nonequispaced samples and Chicken (2007) for the regression model with random uniform design.

This paper presents an extension of a result established by Chicken (2007). We prove that a generalized version of the BlockShrink construction achieves better rates of convergence than the conventional term-by-term thresholding estimators.

The main contributions of this study concern the two following points :

- *The model:* we consider the regression model with (known) random design, not necessarily uniform.
- *The statistical approach:* we adopt the minimax approach under the $\mathbb{L}^p$ risk over Besov balls (regular, sparse and critical zones). The parameter $p$ can be greater than or equal to 2.

From a technical point of view, the proof is significantly more complicated than for the uniform design and the case $p = 2$. We combine a general theorem established by Chesneau (2006, Theorem 4.2) with several probability inequalities such that the Talagrand inequality and the Borel inequality.

The paper is organized as follows. Section 2 introduces the model, the adopted minimax approach, the wavelet bases and the considered estimator. Section 3 presents the main result while Section 4 contains a detailed proof of the main result.

## 2. Model, Wavelet bases and Estimator

### 2.1. The model

We observe $n$ pairs of random variables $\{(X_1, Y_1), \ldots, (X_n, Y_n)\}$ governed by the equation:

$$Y_i = f(X_i) + z_i, \qquad i = 1, \ldots, n, \qquad (2.1)$$

where the design variables $(X_1, \ldots, X_n)$ are i.i.d. with $X_i \in [0, 1]$, the variables $(z_1, \ldots, z_n)$ are i.i.d. Gaussian with mean zero, variance one and are independent of $(X_1, \ldots, X_n)$. We denote by $g$ the density of $X_1$. The function $f$ is unknown. The goal is to estimate $f$ from the observations $\{(X_1, Y_1), \ldots, (X_n, Y_n)\}$. Additional assumptions on the functions $f$ and $g$ will be specified latter (see Theorem 3.1 below).

To estimate $f$, several adaptive methods have been elaborated (according to the nature of the design). See, for instance, the transformation method of Cai and Brown (1998), the model selection method of Baraud (2002) and the kernel method of Gaïffas (2006).

In this study, we shall consider a particular wavelet thresholding estimator. For the sake of clarity, let us denote this estimator by $\hat{f}_n$. The performances of $\hat{f}_n$ will be measured under the global $\mathbb{L}^p$ risk defined by

$$R(\hat{f}_n, f) = \mathbb{E}_f^n(\|\hat{f}_n - f\|_p^p) = \mathbb{E}_f^n \left( \int_0^1 |\hat{f}_n(t) - f(t)|^p dt \right).$$

Here, $p$ is a real number greater than or equal to 2 and $\mathbb{E}_f^n$ is the expectation with respect to the distribution of the observations $\{(X_1, Y_1), \ldots, (X_n, Y_n)\}$. The unknown regression function $f$ is supposed to belong to a wide class of functions: the Besov balls. Wavelets and Besov balls are presented in the next subsection.

### 2.2. Wavelets and Besov balls

We consider an orthonormal wavelet basis generated by dilation and translation of a compactly supported "father" wavelet $\phi$ and a compactly supported "mother" wavelet $\psi$. For the purpose of this paper, we use the periodized wavelet bases on the unit interval. Let us set

$$\phi_{j,k}(x) = 2^{j/2}\phi(2^j x - k), \qquad \psi_{j,k}(x) = 2^{j/2}\psi(2^j x - k).$$

And let us denote the periodized wavelets by $\phi_{j,k}^{per}(x) = \sum_{l \in \mathbb{Z}} \phi_{j,k}(x - l)$, $\psi_{j,k}^{per}(x) = \sum_{l \in \mathbb{Z}} \psi_{j,k}(x - l)$, $x \in [0,1]$. Then, there exists an integer $\tau$ such that the collection $\zeta = \{\phi_{\tau,k}^{per}(x), k = 0, \ldots, 2^\tau - 1; \psi_{j,k}^{per}(x), j = \tau, \ldots, \infty, k = 0, \ldots, 2^j - 1\}$ constitutes an orthonormal basis of $\mathbb{L}^2([0,1])$. The superscript "*per*" will be suppressed from the notations for convenience.

For any integer $l \geq \tau$, a square-integrable function on $[0,1]$ can be expanded into a wavelet series

$$f(x) = \sum_{k=0}^{2^l-1} \alpha_{l,k} \phi_{l,k}(x) + \sum_{j=l}^{\infty} \sum_{k=0}^{2^j-1} \beta_{j,k} \psi_{j,k}(x),$$

where $\alpha_{j,k} = \int_0^1 f(x)\phi_{j,k}(x)dx$ and $\beta_{j,k} = \int_0^1 f(x)\psi_{j,k}(x)dx$.
For further details about wavelets, see Meyer (1990) and Cohen et al. (1993).

Since $\psi$ is compactly supported, the following property of concentration holds: for any $m > 0$, any $j \geq \tau$ and any $x \in [0,1]$, there exists a constant $C > 0$ satisfying

$$\sum_{k=0}^{2^j-1} |\psi_{j,k}(x)|^m \leq C 2^{jm/2}. \tag{2.2}$$

Now, let us define the main sets of function considered in our statistical approach. Let $M \in (0, \infty)$, $s \in (0, \infty)$, $\pi \in [1, \infty]$ and $r \in [1, \infty]$. Let us set $\beta_{\tau-1,k} = \alpha_{\tau,k}$. We say that a function $f$ belongs to the Besov balls $B_{\pi,r}^s(M)$ if and only if there exists a constant $M^* > 0$ such that the associated wavelet coefficients satisfy

$$\left[\sum_{j=\tau-1}^{\infty} [2^{j(s+1/2-1/\pi)}[\sum_{k=0}^{2^j-1}|\beta_{j,k}|^\pi]^{1/\pi}]^r\right]^{1/r} \leq M^*,$$

with the usual modification if $q = \infty$. We work with the Besov balls because of their exceptional expressive power. For a particular choice of parameters $s$, $\pi$ and $r$, they contain the Hölder and Sobolev balls (see, for instance, Meyer (1990)).

### 2.3. The estimator

We are now in the position to describe the main estimator of the study. It is a $\mathbb{L}^p$ version of the BlockShrink estimator developed by Cai (2002) for the Gaussian sequence model.

Let $p \in [2, \infty)$, $d \in (0, \infty)$ and $L$ be the integer $L = \lfloor (\log n)^{p/2} \rfloor$ where the square brackets denote the floor function. Let $j_1$ and $j_2$ be the integers defined by
$$j_1 = \lfloor (p/2) \log_2(\log n) \rfloor, \qquad j_2 = \lfloor (1/2) \log_2(n/\log n) \rfloor.$$
For any $j \in \{j_1, \ldots, j_2\}$, let us set $\mathcal{A}_j = \{1, \ldots, 2^j L^{-1}\}$ and, for any $K \in \mathcal{A}_j$,
$$\mathcal{B}_{j,K} = \{k \in \{0, \ldots, 2^j - 1\} : \ (K-1)L \le k \le KL - 1\}.$$

We define the ($\mathbb{L}^p$ version of the) BlockShrink estimator by

$$\hat{f}_n(x) = \sum_{k=0}^{2^{j_1}-1} \hat{\alpha}_{j_1,k} \phi_{j_1,k}(x) + \sum_{j=j_1}^{j_2} \sum_{K \in \mathcal{A}_j} \sum_{k \in \mathcal{B}_{j,K}} \hat{\beta}_{j,k} \mathbf{1}_{\{\hat{b}_{j,K} \ge dn^{-1/2}\}} \psi_{j,k}(x), \quad (2.3)$$

where $\hat{b}_{j,K} = \left[ L^{-1} \sum_{k \in \mathcal{B}_{j,K}} |\hat{\beta}_{j,k}|^p \right]^{1/p}$ and

$$\hat{\alpha}_{j,k} = n^{-1} \sum_{i=1}^{n} Y_i g(X_i)^{-1} \phi_{j,k}(X_i), \qquad \hat{\beta}_{j,k} = n^{-1} \sum_{i=1}^{n} Y_i g(X_i)^{-1} \psi_{j,k}(X_i). \ (2.4)$$

This estimator was first defined in this $\mathbb{L}^p$ form by Picard and Tribouley (2000) for general statistical models.

*Comments:*

- The sets $\mathcal{A}_j$ and $\mathcal{B}_{j,K}$ are chosen such that $\cup_{K \in \mathcal{A}_j} \mathcal{B}_{j,K} = \{0, \ldots, 2^j - 1\}$, $\mathcal{B}_{j,K} \cap \mathcal{B}_{j,K'} = \emptyset$ for any $K \ne K'$ with $K, K' \in \mathcal{A}_j$, and $|\mathcal{B}_{j,K}| = L = \lfloor (\log n)^{p/2} \rfloor$.
- It is easy to show that $\hat{\alpha}_{j,k}$ and $\hat{\beta}_{j,k}$ are unbiased estimators of $\alpha_{j,k}$ and $\beta_{j,k}$, the wavelet coefficients of $f$. Moreover, they satisfy several probability inequalities which will be at the heart of the proof of the main result. Further details concerning these inequalities are given in Section 4.
- The considered BlockShrink estimator is adaptive since it does not depend on the smoothness of the unknown function $f$. However, it depends on the norm parameter $p$. An open question is : can we construct a block thresholding procedure that is adaptive to the $\mathbb{L}^p$ risk for all p ?

### 3. Main result

Theorem 3.1 below determines the rates of convergence achieved by the Block-Shrink estimator under the $\mathbb{L}^p$ risk over Besov balls.

**Theorem 3.1.** *Let us consider the regression model with random design (2.1). Suppose that :*

- *the unknown regression function $f$ is bounded from above, i.e. $\|f\|_\infty \leq M'$, where $M' > 0$ denotes a known constant.*
- *the density $g$ of $X_1$ is known and bounded from above and below.*

*Let us consider the BlockShrink estimator $\hat{f}_n$ defined by (2.3) with a large enough threshold constant $d$. Let $p \in [2, \infty[$. Then there exists a constant $C > 0$ such that, for any $\pi \in [1, \infty]$, $r \in [1, \infty]$, $s \in (1/\pi + 1/2, \infty)$ and $n$ large enough, we have*

$$\sup_{f \in B^s_{\pi,r}(M)} \mathbb{E}^n_f \left( \|\hat{f}_n - f\|^p_p \right) \leq C\varphi_n,$$

*where*

$$\varphi_n = \begin{cases} n^{-\alpha_1 p}(\log n)^{\alpha_1 p \mathbf{1}_{\{p > \pi\}}}, & \text{when} \quad \epsilon > 0, \\ (\log n/n)^{\alpha_2 p}(\log n)^{(p - \pi/r) + 1_{\{\epsilon = 0\}}}, & \text{when} \quad \epsilon \leq 0, \end{cases}$$

*with $\alpha_1 = s/(2s+1)$, $\alpha_2 = (s - 1/\pi + 1/p)/(2(s - 1/\pi) + 1)$ and $\epsilon = \pi s + 2^{-1}(\pi - p)$.*

The rates of convergence presented in Theorem 3.1 above are minimax except in the cases $\{p > \pi\} \cap \{\epsilon > 0\}$ and $\epsilon = 0$ where there is an extra logarithmic term. They are better than those achieved by the conventional term-by-term thresholding estimators (hard, soft,... ). The main difference is for the case $\{\pi \geq p\}$ where there is no extra logarithmic term. Let us mention that Theorem 3.1 can be proved for $p \in (1, 2)$ if we consider the BlockShrink estimator (2.3) defined with $L = \ln n$. Further details can be found in Chesneau (2006). Further details about the rates of convergence for the regression problem (2.1) via the minimax approach under the $\mathbb{L}^p$ risk over Besov balls can be found in Chesneau (2007).

As mentioned in the motivations of the paper, Theorem 3.1 is an extension of a result proved by Chicken (2007, Theorem 2) for the uniform design, the $\mathbb{L}^2$ risk and the Hölder balls $B^s_{\infty,\infty}(M)$.

*Comments on the choice of the thresholding constant $d$.* From a theoretical point of view, it is difficult to determine the exact minimum value of $d$ such that $\hat{f}_n$ achieves the rates of convergence exhibited in Theorem 3.1. In fact, Theorem 3.1 holds for $d \geq \mu_1$ where $\mu_1$ refers to the constant of Proposition 4.1 below.

## 4. Proof of Theorem 3.1

Thanks to a result proved by Chesneau (2006, Theorem 4.2), the proof of Theorem 3.1 is an immediate consequence of Proposition 4.1 below. This proposition shows that the estimators $(\hat{\beta}_{j,k})_k$ defined by (2.4) satisfy a standard moments inequality and a specific concentration inequality.

**Proposition 4.1.** *Let $p \geq 2$. There exist two constants $\mu_1 > 0$ and $C > 0$ such that, for any $j \in \{j_1, \ldots, j_2\}$, $K \in \mathcal{A}_j$ and $n$ large enough, the estimators $(\hat{\beta}_{j,k})_k$ defined by (2.4) satisfy*

- *the following moments condition:*

$$\mathbb{E}^n_f \left( |\hat{\beta}_{j,k} - \beta_{j,k}|^{2p} \right) \leq Cn^{-p}, \tag{4.1}$$

- *the following concentration inequality:*

$$\mathbb{P}^n_f \left( [L^{-1} \sum_{k \in \mathcal{B}_{j,K}} |\hat{\beta}_{j,k} - \beta_{j,k}|^p]^{1/p} \geq \mu_1 2^{-1} n^{-1/2} \right) \leq 4n^{-p}. \tag{4.2}$$

The moments inequality has been proved by Chesneau (2007). The proof of the concentration inequality (4.2) combines several concentration inequalities such that the Talagrand inequality and the Borel inequality. They are recalled in the two auxiliary lemmas below.

**Lemma 4.1** (Talagrand (1994)). *Let $(V_1, \ldots, V_n)$ be i.i.d. random variables and $(\epsilon_1, \ldots, \epsilon_n)$ be independent Rademacher variables, also independent of $(V_1, \ldots, V_n)$. Let $\mathcal{F}$ be a class of functions uniformly bounded by $T$. Let $r_n : \mathcal{F} \to \mathbb{R}$ be the operator defined by:*

$$r_n(h) = n^{-1} \sum_{i=1}^{n} h(V_i) - \mathbb{E}(h(V_1)).$$

*Suppose that*

$$\sup_{h \in \mathcal{F}} Var(h(V_1)) \leq v, \qquad \mathbb{E} \left( \sup_{h \in \mathcal{F}} \sum_{i=1}^{n} \epsilon_i h(V_i) \right) \leq nH.$$

*Then, there exist two absolute constants $C_1 > 0$ and $C_2 > 0$ such that, for any $t > 0$, we have:*

$$\mathbb{P} \left( \sup_{h \in \mathcal{F}} r_n(h) \geq t + C_2 H \right) \leq \exp \left( -nC_1 \left( t^2 v^{-1} \wedge tT^{-1} \right) \right).$$

**Lemma 4.2** (The Borel inequality (see Adler (1990))). *Let $\mathcal{D}$ be a subset of $\mathbb{R}$. Let $(\eta_t)_{t \in \mathcal{D}}$ be a centered Gaussian process. Suppose that*

$$\mathbb{E} \left( \sup_{t \in \mathcal{D}} \eta_t \right) \leq N, \qquad \sup_{t \in \mathcal{D}} Var(\eta_t) \leq Q.$$

*Then, for any $x > 0$, we have*

$$\mathbb{P} \left( \sup_{t \in \mathcal{D}} \eta_t \geq x + N \right) \leq \exp \left( -x^2/(2Q) \right). \tag{4.3}$$

We are now in the position to prove Proposition 4.1. Here, $C$ represents a constant which may be different from one term to the other. We suppose that $n$ is large enough.

*Proof of the Proposition 4.1.* By the definiton of $\hat{\beta}_{j,k}$, we have the following decomposition
$$\hat{\beta}_{j,k} - \beta_{j,k} = A_{j,k} + B_{j,k},$$
where
$$A_{j,k} = n^{-1}\sum_{i=1}^{n} f(X_i)g(X_i)^{-1}\psi_{j,k}(X_1) - \mathbb{E}_f^n\left(f(X_1)g(X_1)^{-1}\psi_{j,k}(X_1)\right), \quad (4.4)$$

$$B_{j,k} = n^{-1}\sum_{i=1}^{n} g(X_i)^{-1}\psi_{j,k}(X_i)z_i. \quad (4.5)$$

By the $l_p$-Minkowski inequality, for any $\mu > 0$, we have

$$\mathbb{P}_f^n\left([L^{-1}\sum_{k\in\mathcal{B}_{j,K}} |\hat{\beta}_{j,k} - \beta_{j,k}|^p]^{1/p} \geq 2^{-1}\mu n^{-1/2}\right) \leq \mathcal{U} + \mathcal{V},$$

where
$$\mathcal{U} = \mathbb{P}_f^n\left([L^{-1}\sum_{k\in\mathcal{B}_{j,K}} |A_{j,k}|^p]^{1/p} \geq 4^{-1}\mu n^{-1/2}\right),$$

$$\mathcal{V} = \mathbb{P}_f^n\left([L^{-1}\sum_{k\in\mathcal{B}_{j,K}} |B_{j,k}|^p]^{1/p} \geq 4^{-1}\mu n^{-1/2}\right).$$

Let us investigate separately the upper bounds of $\mathcal{U}$ and $\mathcal{V}$.

• *The upper bound for $\mathcal{U}$.* Our goal is to apply the Talagrand inequality described in Lemma 4.1. Let us consider the set $\mathcal{C}_q$ defined by

$$\mathcal{C}_q = \left\{a = (a_{j,k}) \in \mathbb{Z}^*; \ \sum_{k\in\mathcal{B}_{j,K}} |a_{j,k}|^q \leq 1\right\}$$

and the functions class $\mathcal{F}$ defined by

$$\mathcal{F} = \left\{h; \ h(x) = f(x)g(x)^{-1}\sum_{k\in\mathcal{B}_{j,K}} a_{j,k}\psi_{j,k}(x), \ a \in \mathcal{C}_q\right\}.$$

By an argument of duality, we have

$$[\sum_{k\in\mathcal{B}_{j,K}} |A_{j,k}|^p]^{1/p} = \sup_{a\in\mathcal{C}_q} \sum_{k\in\mathcal{B}_{j,K}} a_{j,k}A_{j,k} = \sup_{h\in\mathcal{F}} r_n(h),$$

where $A_{j,k}$ is defined by (4.4) and $r_n$ denotes the function defined in Lemma 4.1. Now, let us evaluate the parameters $T$, $H$ and $v$ of the Talagrand inequality.

First of all, notice that, for $p \geq 2$ (and, a fortiori, $q = 1 + (p-1)^{-1} \leq 2$), an elementary inequality of $l_p$ norm gives

$$\sup_{a \in \mathcal{C}_q} [\sum_{k \in \mathcal{B}_{j,K}} |a_{j,k}|^2]^{1/2} \leq \sup_{a \in \mathcal{C}_q} [\sum_{k \in \mathcal{B}_{j,K}} |a_{j,k}|^q]^{1/q} \leq 1.$$

– *The value of $T$*. Let $h$ be a function in $\mathcal{F}$. By the Cauchy-Schwarz inequality, the assumptions of boundedness of $f$ and $g$ and the property of concentration (2.2), for any $x \in [0,1]$, we find

$$\begin{aligned} |h(x)| &\leq |f(x)||g(x)|^{-1} [\sum_{k \in \mathcal{B}_{j,K}} |\psi_{j,k}(x)|^2]^{1/2} \sup_{a \in \mathcal{C}_q} [\sum_{k \in \mathcal{B}_{j,K}} |a_{j,k}|^2]^{1/2} \\ &\leq \|f\|_\infty \|1/g\|_\infty [\sum_{k \in \mathcal{B}_{j,K}} |\psi_{j,k}(x)|^2]^{1/2} \leq C 2^{j/2}. \end{aligned}$$

Hence $T = C 2^{j/2}$.

– *The value of $H$*. The $l_p$-Hölder inequality and the Hölder inequality imply

$$\begin{aligned} & \mathbb{E}_f^n \left( \sup_{a \in \mathcal{C}_q} \sum_{i=1}^n \sum_{k \in \mathcal{B}_{j,K}} a_{j,k} \epsilon_i f(X_i) g(X_i)^{-1} \psi_{j,k}(X_i) \right) \\ &\leq \sup_{a \in \mathcal{C}_q} [\sum_{k \in \mathcal{B}_{j,K}} |a_{j,k}|^q]^{1/q} \left[ \sum_{k \in \mathcal{B}_{j,K}} \mathbb{E}_f^n \left( |\sum_{i=1}^n \epsilon_i f(X_i) g(X_i)^{-1} \psi_{j,k}(X_i)|^p \right) \right]^{1/p} \\ &\leq \left[ \sum_{k \in \mathcal{B}_{j,K}} \mathbb{E}_f^n \left( |\sum_{i=1}^n \epsilon_i f(X_i) g(X_i)^{-1} \psi_{j,k}(X_i)|^p \right) \right]^{1/p}. \end{aligned} \quad (4.6)$$

Since $(\epsilon_1, \ldots, \epsilon_n)$ are independent Rademacher variables, also independent of $\mathbb{X} = (X_1, \ldots, X_n)$, the Khintchine inequality yields

$$\begin{aligned} & \mathbb{E}_f^n \left( |\sum_{i=1}^n \epsilon_i f(X_i) g(X_i)^{-1} \psi_{j,k}(X_i)|^p \right) \\ &= \mathbb{E}_f^n \left( \mathbb{E}_f^n (|\sum_{i=1}^n \epsilon_i f(X_i) g(X_i)^{-1} \psi_{j,k}(X_i)|^p | \mathbb{X}) \right) \\ &\leq C \mathbb{E}_f^n \left( |\sum_{i=1}^n |f(X_i)|^2 |g(X_i)|^{-2} |\psi_{j,k}(X_i)|^2|^{p/2} \right) = CI. \end{aligned} \quad (4.7)$$

Now, let us consider the i.i.d. random variables $(N_1, \ldots, N_n)$ defined by

$$N_i = |f(X_i)|^2 |g(X_i)|^{-2} |\psi_{j,k}(X_i)|^2, \qquad i = 1, \ldots, n.$$

An elementary inequality of convexity implies $I \leq 2^{p/2-1}(I_1 + I_2)$ where

$$I_1 = \mathbb{E}_f^n\left(|\sum_{i=1}^n (N_i - \mathbb{E}_f^n(N_1))|^{p/2}\right), \quad I_2 = n^{p/2}\mathbb{E}_f^n(N_1)^{p/2}.$$

Let us analyze the upper bounds for $I_1$ and $I_2$, in turn.

– *The upper bound for $I_1$.* The Rosenthal inequality applied to $(N_1, \ldots, N_n)$ and the Cauchy-Schwartz inequality imply

$$\begin{aligned} I_1 &\leq C\left(n\mathbb{E}_f^n(|N_1 - \mathbb{E}_f^n(N_1)|^{p/2}) + \left(n\mathbb{E}_f^n(|N_1 - \mathbb{E}_f^n(N_1)|^2)\right)^{p/4}\right) \\ &\leq C\left(n\mathbb{E}_f^n(|N_1|^{p/2}) + \left(n\mathbb{E}_f^n(|N_1|^2)\right)^{p/4}\right). \end{aligned}$$

For any $m \geq 1$, $j \in \{j_1, \ldots, j_2\}$ and $k \in \{0, \ldots, 2^j - 1\}$, the assumptions of boundedness of $f$ and $g$ give

$$\begin{aligned} \mathbb{E}_f^n(|N_1|^m) &= \int_0^1 |f(x)|^{2m}|g(x)|^{-2m+1}|\psi_{j,k}(x)|^{2m}dx \\ &\leq \|f\|_\infty^{2m}\|\psi\|_\infty^{2m-2}2^{j(m-1)}\int_0^1 |\psi_{j,k}(x)|^2 dx \leq C2^{j_2(m-1)} \leq Cn^{m-1}. \end{aligned}$$

Therefore $I_1 \leq Cn^{p/2}$.

– *The upper bound for $I_2$.* Since $\mathbb{E}_f^n(N_1) \leq C$, we have $I_2 \leq Cn^{p/2}$.

Combining the obtained upper bounds for $I_1$ and $I_2$, we find

$$I \leq C(I_1 + I_2) \leq Cn^{p/2}. \tag{4.8}$$

Putting (4.6), (4.7) and (4.8) together, we see that

$$\mathbb{E}_f^n\left(\sup_{a\in\mathcal{C}_q}\sum_{i=1}^n\sum_{k\in\mathcal{B}_{j,K}} a_{j,k}\epsilon_i f(X_i)g(X_i)^{-1}\psi_{j,k}(X_i)\right) \leq C[\sum_{k\in\mathcal{B}_{j,K}} I]^{1/p} \leq Cn^{1/2}L^{1/p}.$$

Hence $H = Cn^{-1/2}L^{1/p}$.

– *The value of $v$.* By the assumptions of boundedness of $f$ and $g$ and the orthonormality of $\zeta$, we obtain

$$\begin{aligned} \sup_{h\in\mathcal{F}} Var(h(X_1)) \\ &\leq \sup_{a\in\mathcal{C}_q}\mathbb{E}_f^n\left(|f(X_1)|^2|g(X_1)|^{-2}|\sum_{k\in\mathcal{B}_{j,K}} a_{j,k}\psi_{j,k}(X_1)|^2\right) \\ &\leq \|f\|_\infty^2\|1/g\|_\infty \sup_{a\in\mathcal{C}_q}\mathbb{E}_f^n\left(\sum_{k\in\mathcal{B}_{j,K}}\sum_{k'\in\mathcal{B}_{j,K}} a_{j,k}a_{j,k'}g(X_1)^{-1}\psi_{j,k}(X_1)\psi_{j,k'}(X_1)\right) \\ &= C\sup_{a\in\mathcal{C}_q}\sum_{k\in\mathcal{B}_{j,K}}\sum_{k'\in\mathcal{B}_{j,K}} a_{j,k}a_{j,k'}\int_0^1 \psi_{j,k}(x)\psi_{j,k'}(x)dx \\ &= C\sup_{a\in\mathcal{C}_q}\sum_{k\in\mathcal{B}_{j,K}} |a_{j,k}|^2 \leq C. \end{aligned}$$

Hence $v = C$.

Now, let us notice that, for any $j \in \{j_1, \ldots, j_2\}$, we have $n2^j \leq n2^{j_2} \leq 2n^{3/2}(\log n)^{-1/2}$. Since $(\log n)^{1/2} \leq L^{1/p} < 2^{1/p}(\log n)^{1/2}$, for $t = 8^{-1}\mu L^{1/p} n^{-1/2}$, we have

$$\left(t^2 v^{-1} \wedge tT^{-1}\right) \geq C \left(\mu^2(\log n/n) \wedge \mu\sqrt{(\log n/(n2^j))}\right) \geq C\mu^2(\log n/n).$$

So, for $\mu$ large enough and $t = 8^{-1}\mu L^{1/p} n^{-1/2}$, the Talagrand inequality yields

$$\begin{aligned}
\mathcal{U} &= \mathbb{P}_f^n \left([L^{-1} \sum_{k \in \mathcal{B}_{j,K}} |A_{j,k}|^p]^{1/p} \geq 4^{-1}\mu n^{-1/2}\right) \\
&\leq \mathbb{P}_f^n \left([L^{-1} \sum_{k \in \mathcal{B}_{j,K}} |A_{j,k}|^p]^{1/p} \geq 8^{-1}\mu n^{-1/2} + Cn^{-1/2}\right) \\
&\leq \mathbb{P}_f^n \left(\sup_{h \in \mathcal{F}} r_n(h) \geq t + C_2 H\right) \\
&\leq \exp\left(-nC_1\left(t^2 v^{-1} \wedge tT^{-1}\right)\right) \leq \exp\left(-nC\mu^2(\log n/n)\right) \leq n^{-p}.
\end{aligned}$$

We obtain the desired upper bound for $\mathcal{U}$.

• *The upper bound for $\mathcal{V}$.* Here, we apply the Borel inequality described in Lemma 4.2. Let us consider the set $\mathcal{C}_q$ defined by

$$\mathcal{C}_q = \left\{a = (a_{j,k}) \in \mathbb{Z}^*;\ \sum_{k \in \mathcal{B}_{j,K}} |a_{j,k}|^q \leq 1\right\}$$

and the process $\mathcal{Z}(a)$ defined by

$$\mathcal{Z}(a) = \sum_{k \in \mathcal{B}_{j,K}} a_{j,k} B_{j,k},$$

where $B_{j,k}$ is defined by (4.5). Let us notice that, conditionally on $\mathbb{X} = (X_1, \ldots, X_n)$, $\mathcal{Z}(a)$ is a centered Gaussian process. Moreover, by an argument of duality, we have

$$\sup_{a \in \mathcal{C}_q} \mathcal{Z}(a) = \sup_{a \in \mathcal{C}_q} \sum_{k \in \mathcal{B}_{j,K}} a_{j,k} B_{j,k} = [\sum_{k \in \mathcal{B}_{j,K}} |B_{j,k}|^p]^{1/p}.$$

Now, let us investigate separately the upper bounds for $\mathbb{E}_f^n(\sup_{a \in \mathcal{C}_q} \mathcal{Z}(a)|\mathbb{X})$ and $\sup_{a \in \mathcal{C}_q} Var_f^n(\mathcal{Z}(a)|\mathbb{X})$.

– *The upper bound for $\mathbb{E}_f^n(\sup_{a \in \mathcal{C}_q} \mathcal{Z}(a)|\mathbb{X})$.* Let us consider the set $\mathcal{B}_\mu$ defined by

$$\mathcal{B}_\mu = \left\{|n^{-1} \sum_{i=1}^n g(X_i)^{-1}|\psi_{j,k}(X_i)|^2 - 1| \geq \mu\right\}.$$

Let us work on the set $\mathcal{B}_\mu^c$, the complementary of $\mathcal{B}_\mu$. By the Jensen inequality, the fact that $\mathcal{Z}(a) \mid \mathbb{X} \sim \mathcal{N}(0, n^{-2} \sum_{i=1}^n |g(X_i)|^{-2}|\psi_{j,k}(X_i)|^2)$ and the assumptions of boundedness made on $g$, we find

$$\begin{aligned}
\mathbb{E}_f^n\left(\sup_{a\in\mathcal{C}_q} \mathcal{Z}(a)|\mathbb{X}\right) &\leq [\sum_{k\in\mathcal{B}_{j,K}} \mathbb{E}_f^n(|B_{j,k}|^p|\mathbb{X})]^{1/p} \\
&= C[\sum_{k\in\mathcal{B}_{j,K}} (n^{-2}\sum_{i=1}^n |g(X_i)|^{-2}|\psi_{j,k}(X_i)|^2)^{p/2}]^{1/p} \\
&\leq C\|1/g\|_\infty[\sum_{k\in\mathcal{B}_{j,K}} (n^{-2}\sum_{i=1}^n g(X_i)^{-1}|\psi_{j,k}(X_i)|^2)^{p/2}]^{1/p} \\
&\leq Cn^{-1/2}[\sum_{k\in\mathcal{B}_{j,K}} (n^{-1}\sum_{i=1}^n g(X_i)^{-1}|\psi_{j,k}(X_i)|^2 - 1 + 1)^{p/2}]^{1/p} \\
&\leq Cn^{-1/2}[\sum_{k\in\mathcal{B}_{j,K}} (\mu+1)^{p/2}]^{1/p} \leq C(\mu+1)^{1/2}L^{1/p}n^{-1/2}.
\end{aligned}$$

Hence $N = N(\mathbb{X}) = C(\mu+1)^{1/2}L^{1/p}n^{-1/2}$.

− *The upper bound for* $\sup_{a\in\mathcal{C}_q} Var_f^n(\mathcal{Z}(a)|\mathbb{X})$. Let us define the set $\mathcal{A}_\mu$ by

$$\mathcal{A}_\mu = \left\{\sup_{a\in\mathcal{C}_q}(\sum_{k\in\mathcal{B}_{j,K}}\sum_{k'\in\mathcal{B}_{j,K}} a_{j,k}a_{j,k'}(n^{-1}\sum_{i=1}^n g(X_i)^{-1}\psi_{j,k}(X_i)\psi_{j,k'}(X_i)) - \sum_{k\in\mathcal{B}_{j,K}}|a_{j,k}|^2) \geq \mu\right\}.$$

Let us work on the set $\mathcal{A}_\mu^c$, the complementary of $\mathcal{A}_\mu$. Using the assumptions of boundedness of $g$, we have

$$\begin{aligned}
G &= \sup_{a\in\mathcal{C}_q}\left(\sum_{k\in\mathcal{B}_{j,K}}\sum_{k'\in\mathcal{B}_{j,K}} a_{j,k}a_{j,k'}(n^{-1}\sum_{i=1}^n |g(X_i)|^{-2}\psi_{j,k}(X_i)\psi_{j,k'}(X_i))\right) \\
&\leq C\left[\sup_{a\in\mathcal{C}_q}\left(\sum_{k\in\mathcal{B}_{j,K}}\sum_{k'\in\mathcal{B}_{j,K}} a_{j,k}a_{j,k'}\left(n^{-1}\sum_{i=1}^n g(X_i)^{-1}\psi_{j,k}(X_i)\psi_{j,k'}(X_i)\right)\ldots \right.\right. \\
&\quad \left.\left. - \sum_{k\in\mathcal{B}_{j,K}}|a_{j,k}|^2\right) + \sup_{a\in\mathcal{C}_q}\sum_{k\in\mathcal{B}_{j,K}}|a_{j,k}|^2\right] \leq C(\mu+1).
\end{aligned}$$

Since $\mathbb{E}_f^n(z_i z_{i'}) = 1$ if $i = i'$ and $\mathbb{E}_f^n(z_i z_{i'}) = 0$ otherwise, we have

$$\sup_{a \in \mathcal{C}_q} Var_f^n(\mathcal{Z}(a)|\mathbb{X}) = \sup_{a \in \mathcal{C}_q} \mathbb{E}_f^n(|\mathcal{Z}(a)|^2|\mathbb{X})$$

$$= \sup_{a \in \mathcal{C}_q} \mathbb{E}_f^n \left( \sum_{k \in \mathcal{B}_{j,K}} \sum_{k' \in \mathcal{B}_{j,K}} a_{j,k} a_{j,k'} B_{j,k} B_{j,k'} | \mathbb{X} \right)$$

$$= \sup_{a \in \mathcal{C}_q} \left( n^{-2} \sum_{k \in \mathcal{B}_{j,K}} \sum_{k' \in \mathcal{B}_{j,K}} a_{j,k} a_{j,k'} \sum_{i=1}^n \sum_{i'=1}^n |g(X_i)|^{-2} \psi_{j,k}(X_i) \psi_{j,k'}(X_{i'}) \mathbb{E}_f^n(z_i z_{i'}) \right)$$

$$= n^{-1} \sup_{a \in \mathcal{C}_q} \left( \sum_{k \in \mathcal{B}_{j,K}} \sum_{k' \in \mathcal{B}_{j,K}} a_{j,k} a_{j,k'} (n^{-1} \sum_{i=1}^n |g(X_i)|^{-2} \psi_{j,k}(X_i) \psi_{j,k'}(X_i)) \right)$$

$$= n^{-1} G \leq C n^{-1}(\mu + 1).$$

Hence $Q = Q(\mathbb{X}) = C n^{-1}(\mu + 1)$.

The obtained values of $N$ and $Q$ will allow us to conclude. For any $x > 0$, we have

$$\mathbb{P}_f^n \left( \sup_{a \in \mathcal{C}_q} \mathcal{Z}(a) \geq x + C(1+\mu)^{1/2} L^{1/p} n^{-1/2} \right)$$

$$= \mathbb{E}_f^n \left( \mathbb{P}_f^n (\sup_{a \in \mathcal{C}_q} \mathcal{Z}(a) \geq x + C(1+\mu)^{1/2} L^{1/p} n^{-1/2} | \mathbb{X})(1_{\mathcal{B}_\mu} + 1_{\mathcal{B}_\mu^c}) \right)$$

$$\leq \mathbb{P}_f^n(\mathcal{B}_\mu) + \mathbb{E}_f^n \left( \mathbb{P}_f^n (\sup_{a \in \mathcal{C}_q} \mathcal{Z}(a) \geq x + N(\mathbb{X}) | \mathbb{X}) \right). \quad (4.9)$$

The Borel inequality described in Lemma 4.2 implies

$$\mathbb{E}_f^n \left( \mathbb{P}_f^n (\sup_{a \in \mathcal{C}_q} \mathcal{Z}(a) \geq x + N(\mathbb{X}) | \mathbb{X}) \right) \leq \mathbb{E}_f^n \left( \exp(-x^2/(2Q(\mathbb{X}))) \right). \quad (4.10)$$

Moreover, by definition of $\mathcal{A}_\mu$, we have

$$\mathbb{E}_f^n \left( \exp\left(-x^2/(2Q(\mathbb{X}))\right) \right) = \mathbb{E}_f^n \left( \exp\left(-x^2/(2Q(\mathbb{X}))\right)(1_{\mathcal{A}_\mu} + 1_{\mathcal{A}_\mu^c}) \right)$$

$$\leq \mathbb{P}_f^n(\mathcal{A}_\mu) + \exp\left(-nx^2/(2(\mu+1))\right). \quad (4.11)$$

Putting the inequalities (4.9), (4.10) and (4.11) together, for $x = 8^{-1} \mu L^{1/p} n^{-1/2}$ and $\mu$ large enough, we obtain

$$\mathcal{V} = \mathbb{P}_f^n \left( \sup_{a \in \mathcal{C}_q} \mathcal{Z}(a) \geq 4^{-1} \mu L^{1/p} n^{-1/2} \right)$$

$$\leq \mathbb{P}_f^n \left( \sup_{a \in \mathcal{C}_q} \mathcal{Z}(a) \geq 8^{-1} \mu L^{1/p} n^{-1/2} + C(1+\mu)^{1/2} L^{1/p} n^{-1/2} \right)$$

$$\leq \mathbb{P}_f^n(\mathcal{A}_\mu) + \mathbb{P}_f^n(\mathcal{B}_\mu) + \exp\left(-C \mu^2 L^{2/p}/(\mu+1)\right). \quad (4.12)$$

Lemma 4.3 below provides the upper bounds for $\mathbb{P}_f^n(\mathcal{A}_\mu)$ and $\mathbb{P}_f^n(\mathcal{B}_\mu)$.

**Lemma 4.3.** *For $\mu$ and $n$ large enough, we have*

$$\max\left(\mathbb{P}_f^n(\mathcal{A}_\mu), \mathbb{P}_f^n(\mathcal{B}_\mu)\right) \leq n^{-p}.$$

By the inequality (4.12), the fact that $L = \lfloor (\log n)^{p/2} \rfloor$ and Lemma 4.3, for $\mu$ large enough, we have

$$\mathcal{V} \leq 3n^{-p}.$$

Combining the obtained upper bounds for $\mathcal{U}$ and $\mathcal{V}$, we achieve the proof of Proposition 4.1. □

*Proof of Lemma 4.3..* Let us investigate the upper bounds for $\mathbb{P}_f^n(\mathcal{B}_\mu)$ and $\mathbb{P}_f^n(\mathcal{A}_\mu)$.

- *The upper bound for $\mathbb{P}_f^n(\mathcal{B}_\mu)$.* First of all, notice that the random variables

$$(|\psi_{j,k}(X_1)|^2 g(X_1)^{-1}, \ldots, |\psi_{j,k}(X_n)|^2 g(X_n)^{-1}),$$

are i.i.d. and, since $g$ is bounded from below, we have

$$|\psi_{j,k}(X_i)|^2 g(X_i)^{-1} \leq \|1/g\|_\infty \|\psi\|_\infty^2 2^j, \qquad \mathbb{E}_f^n\left(|\psi_{j,k}(X_1)|^2 g(X_1)^{-1}\right) = 1.$$

So, for any $j \in \{j_1, \ldots, j_2\}$, the Hoeffding inequality implies the existence of a constant $C > 0$ such that

$$\mathbb{P}_f^n(\mathcal{B}_\mu) \leq 2\exp\left(-Cn\mu^2 2^{-2j}\right) \leq 2\exp\left(-Cn\mu^2 2^{-2j_2}\right) \leq n^{-C\mu^2}.$$

We obtain the desired upper bound by taking $\mu$ large enough.

- *The upper bound for $\mathbb{P}_f^n(\mathcal{A}_\mu)$.* The goal is to apply the Talagrand inequality described in Lemma 4.1. Let us consider the set $\mathcal{C}_q$ defined by

$$\mathcal{C}_q = \left\{a = (a_{j,k}) \in \mathbb{Z}^*; \sum_{k \in \mathcal{B}_{j,K}} |a_{j,k}|^q \leq 1\right\}$$

and the functions class $\mathcal{F}'$ defined by

$$\mathcal{F}' = \left\{h; \, h(x) = g(x)^{-1} \sum_{k \in \mathcal{B}_{j,K}} \sum_{k' \in \mathcal{B}_{j,K}} a_{j,k} a_{j,k'} \psi_{j,k}(x) \psi_{j,k'}(x), \; a \in \mathcal{C}_q\right\}.$$

We have

$$\sup_{a \in \mathcal{C}_q} \left(\sum_{k \in \mathcal{B}_{j,K}} \sum_{k' \in \mathcal{B}_{j,K}} a_{j,k} a_{j,k'} (n^{-1} \sum_{i=1}^n g(X_i)^{-1} \psi_{j,k}(X_i) \psi_{j,k'}(X_i)) - \sum_{k \in \mathcal{B}_{j,K}} |a_{j,k}|^2\right)$$
$$= \sup_{h \in \mathcal{F}'} r_n(h),$$

where $r_n$ denotes the function defined in Lemma 4.1. Thus, it suffices to determine the parameters $T$, $H$ and $v$ of the Talagrand inequality.

− *The value of T.* Let $h$ be a function of $\mathcal{F}'$. Using the Hölder inequality, the fact that $g$ is bounded from below and the concentration property (2.2), for any $x \in [0,1]$, we find

$$|h(x)| \leq \|1/g\|_\infty \sum_{k \in \mathcal{B}_{j,K}} |a_{j,k}|^2 \sum_{k \in \mathcal{B}_{j,K}} |\psi_{j,k}(x)|^2 \leq C 2^j.$$

Hence $T = C 2^j$.

− *The value of H.* The Cauchy-Schwarz inequality implies that

$$\mathbb{E}_f^n \left( \sup_{a \in \mathcal{C}_q} \sum_{k \in \mathcal{B}_{j,K}} \sum_{k' \in \mathcal{B}_{j,K}} a_{j,k} a_{j,k'} \sum_{i=1}^n \epsilon_i g(X_i)^{-1} \psi_{j,k}(X_i) \psi_{j,k'}(X_i) \right)$$

$$\leq \sup_{a \in \mathcal{C}_q} [\sum_{k \in \mathcal{B}_{j,K}} \sum_{k' \in \mathcal{B}_{j,K}} |a_{j,k}|^2 |a_{j,k'}|^2]^{1/2} \ldots$$

$$\left[ \sum_{k \in \mathcal{B}_{j,K}} \sum_{k' \in \mathcal{B}_{j,K}} \mathbb{E}_f^n \left( |\sum_{i=1}^n \epsilon_i g(X_i)^{-1} \psi_{j,k}(X_i) \psi_{j,k'}(X_i)|^2 \right) \right]^{1/2}$$

$$\leq \left[ \sum_{k \in \mathcal{B}_{j,K}} \sum_{k' \in \mathcal{B}_{j,K}} \mathbb{E}_f^n \left( |\sum_{i=1}^n \epsilon_i g(X_i)^{-1} \psi_{j,k}(X_i) \psi_{j,k'}(X_i)|^2 \right) \right]^{1/2} .(4.13)$$

Since $(\epsilon_1, \ldots, \epsilon_n)$ are independent Rademacher variables, also independent of $\mathbb{X} = (X_1, \ldots, X_n)$, the Khintchine inequality and the fact that $g$ is bounded from below give

$$\mathbb{E}_f^n \left( |\sum_{i=1}^n \epsilon_i g(X_i)^{-1} \psi_{j,k}(X_i) \psi_{j,k'}(X_i)|^2 \right)$$

$$= \mathbb{E}_f^n \left( \mathbb{E}_f^n \left( |\sum_{i=1}^n \epsilon_i g(X_i)^{-1} \psi_{j,k}(X_i) \psi_{j,k'}(X_i)|^2 \right) | \mathbb{X} \right)$$

$$\leq C \mathbb{E}_f^n \left( \sum_{i=1}^n |g(X_i)|^{-2} |\psi_{j,k}(X_i)|^2 |\psi_{j,k'}(X_i)|^2 \right)$$

$$\leq C \|1/g\|_\infty^2 n \mathbb{E}_f^n \left( |\psi_{j,k}(X_1)|^2 |\psi_{j,k'}(X_1)|^2 \right). \tag{4.14}$$

Using the property of concentration (2.2) and the inequalities (4.13)-(4.14), we

find

$$\mathbb{E}_f^n \left( \sup_{a \in \mathcal{C}_q} \sum_{k \in \mathcal{B}_{j,K}} \sum_{k' \in \mathcal{B}_{j,K}} a_{j,k} a_{j,k'} (\sum_{i=1}^n \epsilon_i g(X_i)^{-1} \psi_{j,k}(X_i) \psi_{j,k'}(X_i)) \right)$$

$$\leq C \left[ n \mathbb{E}_f^n \left( |\sum_{k \in \mathcal{B}_{j,K}} |\psi_{j,k}(X_i)|^2|^2 \right) \right]^{1/2} \leq C n^{1/2} 2^j.$$

Hence $H = C 2^j n^{-1/2}$.

− *The value of v.* Using the fact that $g$ is bounded from below, the Hölder inequality and the property of concentration (2.2), we have

$$\sup_{h \in \mathcal{F}} Var(h(X_1)) \leq \sup_{a \in \mathcal{C}_q} \mathbb{E}_f^n \left( |g(X_1)|^{-2} | \sum_{k \in \mathcal{B}_{j,K}} \sum_{k' \in \mathcal{B}_{j,K}} a_{j,k} a_{j,k'} \psi_{j,k}(X_1) \psi_{j,k'}(X_1) |^2 \right)$$

$$\leq C \|1/g\|_\infty^2 \sup_{a \in \mathcal{C}_q} [\sum_{k \in \mathcal{B}_{j,K}} |a_{j,k}|^2]^2 \mathbb{E}_f^n \left( |\sum_{k \in \mathcal{B}_{j,K}} |\psi_{j,k}(X_1)|^2|^2 \right)$$

$$\leq C 2^{2j}.$$

Hence $v = C 2^{2j}$.

Now, let us notice that if $t = 2^{-1} \mu$ then

$$\left( t^2 v^{-1} \wedge t T^{-1} \right) \geq C \left( \mu^2 2^{-2j} \wedge \mu 2^{-j} \right) = C \mu^2 2^{-2j}.$$

For any $j \in \{j_1, \ldots, j_2\}$, $t = 2^{-1} \mu$ with $\mu$ large enough, the Talagrand inequality gives

$$\mathbb{P}_f^n(\mathcal{A}_\mu)$$

$$\leq \mathbb{P}_f^n (\sup_{a \in \mathcal{C}_q} ( \sum_{k \in \mathcal{B}_{j,K}} \sum_{k' \in \mathcal{B}_{j,K}} a_{j,k} a_{j,k'} (n^{-1} \sum_{i=1}^n g(X_i)^{-1} \psi_{j,k}(X_i) \psi_{j,k'}(X_i)) - \ldots$$

$$\sum_{k \in \mathcal{B}_{j,K}} |a_{j,k}|^2) \geq 2^{-1} \mu + C 2^j n^{-1/2}) \leq \mathbb{P}(\sup_{h \in \mathcal{F}} r_n(h) \geq t + C_2 H)$$

$$\leq \exp\left( -n C_1 \left( t^2 v^{-1} \wedge t T^{-1} \right) \right) \leq \exp\left( -n C \mu^2 2^{-2j} \right)$$

$$\leq \exp\left( -n C \mu^2 2^{-2j_2} \right) \leq n^{-p}.$$

This ends the proof of Lemma 4.3. □

**Acknowledgment.** We thank two referees and the associate editor for their thorough and useful comments which have helped to improve the presentation of the paper.